\documentclass[11pt]{article}   
\usepackage{amssymb,amscd,latexsym}   
\usepackage{amsmath}
\usepackage{amsthm}
\newcommand{\rar}{\rightarrow}
\newcommand{\lar}{\longrightarrow}
\newcommand{\llar}{-\kern-5pt-\kern-5pt\longrightarrow}
\newcommand{\surjects}{\twoheadrightarrow}
\newcommand{\injects}{\hookrightarrow}

\newtheorem{Theorem}{Theorem}[section]
\newtheorem{Lemma}[Theorem]{Lemma}
\newtheorem{Corollary}[Theorem]{Corollary}
\newtheorem{Proposition}[Theorem]{Proposition}
\newtheorem{Remark}[Theorem]{Remark}
\newtheorem{Example}[Theorem]{Example}

\newtheorem{Definition}[Theorem]{Definition}
\newtheorem{Question}[Theorem]{Question}
\def\sqr#1#2{{\vcenter{\hrule height.#2pt
        \hbox{\vrule width.#2pt height#1pt \kern#1pt
            \vrule width.#2pt}
        \hrule height.#2pt}}}
\def\phi{\varphi}
\def\demo{\noindent{\bf Proof. }}
\def\square{\mathchoice\sqr64\sqr64\sqr{4}3\sqr{3}3}
\def\qed{\hspace*{\fill} $\square$}
\def\QED{\hbox{\qed}}

\def\xx{{\bf x}}
\def\yy{{\bf y}}

\def\XX{{\bf X}}
\def\YY{{\bf Y}}

\def\ff{{\bf f}}

\def\ff{{\bf f}}
\def\gg{{\bf g}}
\def\hh{{\bf h}}

\def\fn{{\mathfrak n}}

\def\hht{{\rm ht}\,}

\def\edim{{\rm edim}\,}

\def\ker{{\rm ker}\,}
\def\grade{{\rm grade}\,}
\def\rk{{\rm rank}\,}

\def\restr{{\kern-1pt\restriction\kern-1pt}}

\def\NN{\mathbb N}

\def\pp{{\mathbb P}}

\begin{document}
\begin{center}
{\Large{\bf\sc A characteristic free criterion of birationality}}\footnotetext{AMS 1980 Mathematics
Subject Classification (1985   Revision). Primary 14E05, 14E07, 14M05,
13H10; Secondary 14M12, 14M15, 13B10, 13C14, 13F45 .}

\vspace{0.3in}

{\large\sc A. V. D\'oria}\footnote[1]{On leave from UFS, Brazil},
{\large\sc S. H. Hassanzadeh}\footnote[2]{Pos-doctoral Fellow under a CNPq Grant, Brazil}
 and {\large\sc A.  Simis}\footnote[3]{Partially
supported by a CNPq grant, Brazil.}

\end{center}


\bigskip

\begin{abstract}

One develops {\em ab initio} the theory of rational/birational maps
over reduced, but not necessarily irreducible, projective varieties in
arbitrary characteristic. A numerical invariant of a rational map is
introduced,  called the Jacobian dual rank. It is proved that a rational map
in this general setup is birational if and only if the Jacobian dual rank attains
its maximal possible value.
Even in the ``classical'' case where the source variety is irreducible
there is some gain for this invariant over the degree of the map as it
is, on one hand, intrinsically related to natural constructions in commutative algebra
and, on the other hand, is effectively straightforwardly computable.
Applications are given to results so far only known in characteristic zero.
In particular, the surprising result of Dolgachev concerning the degree
of a plane polar Cremona map is given an alternative conceptual angle.

\end{abstract}

\section*{Introduction}

A rational map $\mathfrak{F}:\pp^n\dasharrow \pp^m$ between projective
spaces is defined by $m+1$ forms $\ff=\{f_0,\ldots f_m\}$
of the same degree in $n+1$ variables not all vanishing.
The idea of looking at the syzygies of the forms $\ff$ goes
back to the work of Hulek-Katz-Schreyer (\cite{HKS}) in the case
where $m=n$, as a method of obtaining sufficient conditions for
$\mathfrak{F}$ to be a Cremona map. In \cite{cremona}  this method has been developed
to a greater extent when $m\geq n$. Still more
recently, the present third author pushed the method further to the study
of general rational maps between two integral projective schemes in
arbitrary characteristic by an extended ideal-theoretic method
emphasizing the role of the Rees algebra associated to the ideal
generated by $\ff$ (\cite{bir2003}).

Although the general results were characteristic-free, some
particular questions remained a challenge in  positive
characteristic (but see \cite{SimisVilla} for the case of rational
maps defined by monomials). Even in characteristic
zero a few issues were left answered in these papers.

The goal of this work is both a generalization and a great simplification of the
basic theory of rational maps in arbitrary characteristic.
We consider rational maps defined on reduced, but not necessarily, irreducible
projective varieties.
Though usually thought as a trivial extension from the case where the source
is integral, to establish it on firm ground requires non-trivial proofs on several
passages. Moreover, a {\em bona fide} generalization ought to work well through the
various versions of the notion of birationality. In order to proceed in a rigorous way we
found it essential to define everything {\em ab initio}, so to avoid often inevitable
hand waving from the integrality habits.

The simplification comes about by showing that birationality is
controlled by the behavior of a unique numerical invariant -- called the {\em Jacobian
dual rank} of a rational map.
Alas, this sounds like old knowledge because the classical theory also depends only
on the degree of the rational map. However, the latter is in full control only
in the integral case. Habitually, in positive characteristic one treats the inseparability
degree apart from the main stream of the natural ideas in birational theory.
The new invariant introduced here looks more intrinsic
and makes no explicit reference to inseparability, so the criterion itself and the
applications will be characteristic-free.
Finally, the Jacobian dual rank is straightforwardly effectively computable in the usual implementation of
the Gr\"obner basis algorithm, an appreciable advantage over the field degree.

In addition, the Jacobian dual rank calls attention to several aspects of the  theory
of Rees algebras and base ideals of maps, a trend sufficiently shown in many modern accounts
(see, e.g., \cite{CRS},  \cite{syl2}, \cite{HKS}, \cite{PanRusso}).

\medskip

The paper is divided in two sections.

The first section hinges on the needed background to state the general criterion of birationality.

In the initial subsections we develop the ground material on rational and birational  maps on a reduced source.
Our approach is entirely algebraic, but we mention the transcription to the geometric
side. A degree of care is required to show that the present notion is stable under the expected
manipulations from the ``classical'' case.
One main result in this part is Proposition~\ref{birational_in_terms_of_fields} which drives
us back to an analogue of the field extension version.

The main core is the subsequent subsection, where we introduce the Jacobian dual rank
and prove the basic characteristic-free criterion of birationality  in terms of this rank.
The criterion also holds component-wise as possibly predictable (but not obviously proved!).
We took pains to transcribe the criterion into purely geometric terms, except for the
Jacobian dual rank itself, whose geometric meaning is not entirely apparent at this stage.
This concept has evolved continuously from previous notions, of
which the cradle is the Jacobian dual fibration of \cite{dual} (see
the subsequent \cite{cremona} and \cite{bir2003}).
The base ideal of a rational map will play an essential role often not obviously
foreseen in the basic definitions.

Among the techniques for this part, a giant role  is played by the Rees algebra of an ideal in a more vivid
context.
Thus, let $\ff=\{f_0,\ldots,f_m\}\subset R$ be an ordered set of forms of the same degree,
where $R$ denotes a standard graded $k$-algebra, with $k$ a field and $\dim R\geq 1$.
Since $\ff$ is generated in same degree, the Rees algebra ${\cal R}_R((\ff))$ is a bigraded
standard $k$-algebra. In particular, there is a presentation $R[\YY]=R[Y_0,\ldots,Y_m]\surjects {\cal R}_R((\ff))$
mapping $Y_j\mapsto f_j$.
Let  ${\cal J}=\ker(R[\YY]\surjects {\cal R}_R((\ff))$ denote the corresponding kernel (often
called a presentation ideal of ${\cal R}_R(\ff)$).
It will be seen, among the preliminaries, that $\mathcal{J}$ depends only on the rational map defined by $\ff$
and not on this particular representative.
Note that it is a bihomogeneous ideal
in the natural standard bigrading of $R[\YY]$ induced by the standard
bivariate polynomial ring $k[\XX,\YY]$. Thus, one has
$${\cal J}=\bigoplus_{(p,q)\in \mathbb{N}^2} {\cal J}_{(p,q)},$$
where ${\cal J}_{(p,q)}$ denotes the $k$-vector space of forms of bidegree $(p,q)$.

One first relevant bihomogeneous piece for our purpose is ${\cal J}_{0,*}$, which is
spanned by the forms of bidegree  $(0,q)$, for all $q\geq 1$.
Since these forms have coefficients only in $k$, one can view them as elements of $k[\YY]$.
As such they generate the ideal of all polynomial relations (over $k$) of $\ff$.
Thus,
$$k[\YY]/({\cal J}_{0,*})k[\YY]\simeq k[\ff],$$
as standard graded $k$-algebras,
after a degree renormalization in $k[\ff]$.

For birationality, yet another bihomogeneous piece is as important:
$${\cal J}_{1,*}:=\bigoplus_{r\in\NN} {\cal J}_{1,q}$$
with ${\cal J}_{1,q}$ denoting the bigraded piece of ${\cal J}$ spanned by the forms of bidegree
 $(1,q)$ for all $q\geq 0$. Note that ${\cal J}_{1,*}=\{0\}$
 may very well be the case in general.
 The definition of the Jacobian dual rank is drawn from this datum in form of the rank
 of a suitable matrix. Due to its technicality we refer to the appropriate place in the section.

While elimination theory deals properly with $({\cal J}_{\,0,*})k[\YY]$, our algebraic
approach to birational theory deals with the
puzzle as to how these two  pieces of ${\cal J}$  relate to each other.
Thus, quite a bit of deep classical theory seems to be embedded into these two initial
graded pieces of the defining equations of the blowing up of the base locus.
The purpose of the Jacobian dual rank is roughly to link between these pieces.
While it is not entirely clear how they entangle, the application in this work surfaces
when this rank acquires maximal possible value.

This is the basic outline of the first section.
The main results are contained in Proposition~\ref{rank_and_dimension}
and Theorem~\ref{main_criterion_algebraic}.

\medskip

The second section is concerned with applications of the main criterion.

Since the criterion is stated in terms of an effective numerical invariant,
the obvious applicability is to finding out whether a rational map given explicitly
by a representative is birational and, in the affirmative case, to exhibiting
a representative of the inverse map.

This is not what one has in mind in the section.
Rather, we collected a few results previously known in characteristic zero,
now fully available in arbitrary characteristic.
The flagship is Theorem~\ref{linearmaximalrank}.
One of its consequences is a result for plane Cremona maps, saying that
the only such maps whose base ideal is locally a complete intersection
at its minimal primes are the classical quadratic maps (up to a projective change
of coordinates).
In particular, the result of Dolgachev (\cite{Dol}) is equivalent to showing
that the gradient ideal of a homogeneous homaloidal polynomial in $3$ variables
is locally a complete intersection
at its minimal primes.

We close the section by introducing yet another numerical invariant, called the
semilinear generation degree. It would be interesting to capture its full meaning
in the present and related theories -- such as elimination.

\section{A criterion of birationality}\label{section_criterion}

\subsection{Prolegomenon: the Koszul--Hilbert lemma}

The following curious piece of linear algebra, for which we have found no explicit reference in the literature, will
be useful in the subsequent parts of this section.

\begin{Proposition}\label{koszul_hilbert}
Let $A$ be a commutative ring and $n$ a positive integer.
Let $\rho$ stand for an  $(n+1)\times n$ matrix with arbitrary entries in $A$ and let ${\bf \Delta}=
\{\Delta_0,\ldots,\Delta_n\}$ stand for its {\rm (}ordered, signed{\rm )} $n$-minors.
Then the Koszul relations of ${\bf \Delta}$ belong to the submodule
$${\rm Im}(\rho)\subset A^{n+1}.$$
\end{Proposition}
\demo
Let $B=A[\XX]$ be a polynomial ring in $(n+1)\times n$ indeterminates and let $\mathcal{X}$ denote the corresponding
$(n+1)\times n$ matrix whose entries are $\XX$.

\medskip

We need:

\begin{Lemma}
Let $D$ denote a Noetherian factorial subdomain of $A$ -- e.g., the prime subring of $A$.
Then the ideal of $D[\XX]$ generated by the $n$-minors of $\mathcal{X}$ has grade $\geq 2$.
\end{Lemma}
\demo This is clear as, say, $\Delta_1\in D[\XX]$ generates a prime ideal and, say, $\Delta_2$
is not a multiple of $\Delta_1$.
\qed

\medskip

As a consequence, the well-known Hilbert--Burch principle says that the sequence
$$0\rar C^n\stackrel{\mathcal{X}}{\lar} C^{n+1}\stackrel{{\bf \Delta}}{\lar} C $$
is exact. Clearly, then the Koszul relations of ${\bf \Delta}$ belong to the image of the map
$\mathcal{X}$ as this is exactly the matrix of first syzygies of ${\bf \Delta}$.

Now specialize the variables $\XX$ to the entries of $\rho$ by means of the obvious homomorphism $D[\XX]\lar A$
extending the inclusion $D\subset A$.
Since the formation of determinants and the Koszul complex specialize, it is clear that the specialization preserves
the above relation.
\qed

\begin{Remark}\rm Note that one is not asserting that
$\rho$ is the entire matrix of syzygies of ${\bf \Delta}$, which is only true if the grade
of ${\bf \Delta}$ is at least $2$ as well.
\end{Remark}

\subsection{Preliminaries on rational maps on reduced varieties}

We first fix our standing setup and notation.

Let $k$ denote a fixed, but otherwise arbitrary, ground field.
As usual, in order that a projective set of zeros have enough points one allows coordinates to lie on
an algebraically closed field containing $k$.

Let $R=k[\xx]=k[x_0,\ldots, x_n]\simeq k[\XX]/{\mathfrak a}=k[X_0,\ldots, X_n]/{\mathfrak a}$ denote a standard
graded $k$-algebra over a field $k$ of positive dimension.
In the sequel we always write $\XX$ and $\xx$ for $X_0,\ldots, X_n$ and $x_0,\ldots, x_n$, respectively.
The same notation principle will apply to other sets of variables as needed.

\subsubsection*{Rational data and rational maps}

Let $K(R)$ denote the homogeneous total ring of quotients of $R$ --
i.e., the ring of fractions $\mathfrak{M}^{-1}R$, where $\mathfrak{M}$ is the set of homogeneous regular elements of $R$
not belonging to any of the minimal primes of $R$.

We will be dealing with ordered $(m+1)$-sets $\ff=\{f_0,\ldots, f_m\}$ of forms of the same
degree in $R$.
Let us say that two such ordered $(m+1)$-sets $\ff=\{f_0,\ldots, f_m\}$ and $\ff'=\{f'_0,\ldots, f'_m\}$,
are {\sc equivalent} if the corresponding tuples $(f_0,\ldots ,f_m)$ and $(f'_0,\ldots ,f'_m)$ are proportional
by means of a homogeneous invertible element of $K(R)$.
This means that there exist regular homogeneous elements $F,F'\in R$ such that $Ff'_j=F'f_j$, for $j=0,\ldots,m$.
This easily seen to be an equivalence relation.

Note that equivalent ordered homogeneous $(m+1)$-sets $\ff,\ff'\subset R$ induce a graded $k$-isomorphism of the respective
representative algebras $k[\ff]$ and $k[\ff']$, such that $f_j\longleftrightarrow f'_j$.

\begin{Definition}\label{rational_datum}\rm
A {\sc rational $(m+1)$-datum of degree $d$} on $R$ is an ordered set of $m+1$ forms
$\ff=\{f_0,\ldots f_m\}\subset R$
of the same degree $\geq 1$ satisfying the following conditions:
\begin{enumerate}
\item[{\rm (i)}] $R$ is torsionfree over $k[\ff]$ via the ring extension $k[\ff]\subset R$
\item[{\rm (ii)}] The ideal $(\ff)R\subset R$ is not contained in any minimal prime of $R$.
\end{enumerate}
\end{Definition}
Note that the first condition reads in the present circumstances that every minimal
prime of $R$ contracts to a minimal prime of $k[\ff]$, while the second says that the
ideal $(\ff)\subset R$ has a regular element.

An $m$-{\sc targeted rational map} with {\sc source} Proj$(R)\subset \pp^n_k$ is a collection of equivalent rational $(m+1)$-data.

\begin{Lemma}\label{minors}\rm If $k$ is an infinite field, two homogeneous $(m+1)$-sets $\ff,\ff'\subset R$ satisfying condition (ii) of the above definition are equivalent if and only if the matrix
$$\left(
\begin{array}{cccc}
f_0&f_1&\ldots &f_m\\
f'_0&f'_1&\ldots &f'_m
\end{array}
\right)
$$
has rank $1$.
\end{Lemma}
\demo
Indeed, the rank condition means that the $2\times 2$ minors vanish (which is readily implied
by equivalence as defined) and that the ideal generated by the entries has a regular element (which is obviously the case if one of the ideals generated by the rows has a regular element).
To see the converse, since $k$ is infinite one can choose suitable $\lambda_j\in k$ such that
both $\lambda_0f_0+\cdots +\lambda_mf_m$ and $\lambda_0f'_0+\cdots +\lambda_mf'_m$ are
regular elements.
By hypothesis, $f_jf'_i=f_if'_j$ for all pairs $i,j$.
Therefore, for arbitrary $i$, one has $F'f_i=Ff'_i$, where $F=\lambda_0f_0+\cdots +\lambda_mf_m$
and $F'=\lambda_0f'_0+\cdots +\lambda_mf'_m$.
\qed

\medskip

The rank property was taken in \cite{bir2003} as definition
of equivalent representatives of the same rational map.

\medskip

Any individual rational $(m+1)$-datum is said to be a {\sc representative} of the map and, in an informal
speech, we say that it defines the rational map.
Throughout we denote a rational map by $\mathfrak{F}:\ff$, where $\ff$ is a particular representative.
Often we identify the map with one representative, which will cause no confusion.
As remarked above, equivalent rational $(m+1)$-data $\ff,\ff'\subset R$ induce a $k$-isomorphism of the respective
representative algebras $k[\ff]$ and $k[\ff']$, such that $f_j\longleftrightarrow f'_j$.
Although this isomorphism is graded, it is not homogeneous (of degree $0$). Nevertheless, one can renormalize the gradings of all these algebras
so they become standard $k$-algebras. What this actually means is that one can take
a common polynomial representation $k[\YY]/\mathfrak{b}\simeq k[\ff]$ for all the representatives
of the same $m$-targeted rational map, with $k[\YY]/\mathfrak{b}$ standard.

\subsubsection*{Base ideal and stability}

 The ideal $(\ff)\subset R$ generated by a representative $\ff\subset R$ is called a {\sc base ideal} of the map.
At the level of base ideals, two representatives $\ff$ and $\ff'$ of the same map give an isomorphism of $(\ff)$ and $(\ff')$ as
$R$-submodules of $K(R)$, i.e., as fractional ideals.

The notion of a representative or of a base ideal is pretty stable in the following sense.

\begin{Lemma}\label{stability}
Let $\ff=\{f_0,\ldots f_m\}$ and $\ff'=\{f'_0,\ldots f'_m\}$ be equivalent ordered sets in $R$.
Then $\ff$ is a rational datum if and only if $\ff'$ is a rational datum.
\end{Lemma}
\demo By the obvious symmetry in an argument, it suffices to show one implication.
Thus, assume that $\ff$ satisfies the conditions stated in Definition~\ref{rational_datum}.
As noted above, $k[\ff]\simeq k[\ff']$, so there is a bijection between the respective sets of minimal primes.
This shows that condition (i) of [loc.cit.] holds for $\ff'$.
Since the $k$-isomorphism induces an isomorphism of the ideals $(\ff), (\ff')\subset R$ as $R$-modules (also
remarked before), condition (ii) follows as well.
\qed

\subsubsection*{In terms of total rings of quotients}

 The usual ``field'' criterion reads as well. Namely, choose a fixed presentation $S=k[\YY]/\mathfrak{b}\simeq k[\ff]$
for an arbitrary algebra representative of degree $d$ of an $m$-targeted rational map.
There are actually three versions of the criterion, as follows.
Since $S\hookrightarrow R$ is torsionfree, it induces the following maps:
\begin{itemize}
\item[(1)] An injective graded homomorphism $K(S)\hookrightarrow  K(R)$ of degree $d-1$ of the respective homogeneous total rings
of quotients;
\item[(2)] An injective {\em homogeneous} homomorphism (i.e., a graded homomorphism of degree zero)  $K(S)\hookrightarrow K(R^{(d)})$,
where $R^{(d)}\subset R$ is the $d$th Veronese subring of $R$;
\item[(3)] An injective  (non-graded) homomorphism $K(S)_{_0} \hookrightarrow  K(R)_{_0}$ in degree zero.
\end{itemize}
Note that, by definition, $K(R^{(d)})_{_0}=K(R)_{_0}$, so passing to the Veronese in degree zero is tantamount to doing it ``on the nose''
over $R$.

\subsubsection*{Composition}

 Two such objects can be ``composed'' in the following sense:
\begin{Lemma}\label{composite}
Let $R=k[\XX]/\mathfrak{a}$ be a reduced standard graded $k$-algebra and let $\mathfrak{F}:\ff\subset R$
be a rational $(m+1)$-datum of degree $d$ on $R$.
Take a presentation $S:=k[\YY]/\mathfrak{b}\simeq k[\ff]$, with $\YY=\{y_0,\ldots, y_m\}$ variables over $k$.
Let $\gg:=\{g_0,\ldots, g_r\}\subset S$ denote a rational $(r+1)$-datum of degree $e$ on $S$.
Then:
\begin{enumerate}
\item[{\rm (i)}] $\gg(\ff):=\{g_0(\ff),\ldots,g_r(\ff)\}$ is a rational $(r+1)$-datum of degree $de$ on $R$
\item[{\rm (ii)}] If $\ff'\subset R$ and $\gg'\subset S$ are representatives of $\ff$ and $\gg$, respectively, then
$\gg'(\ff'):=\{g'_0(\ff'),\ldots,g'_r(\ff')\}$ is a representative of the rational map defined by $\gg(\ff)$.
\end{enumerate}
\end{Lemma}
\demo (i) By assumption, $R$ is torsionfree over $k[\ff]\simeq S$ and $S$ is torsionfree over $k[\gg]$.
Therefore, $R$ is torsionfree over $k[\gg]$.
Likewise, by assumption, $(\gg)S$ is not contained in any minimal prime of $S$, i.e., its image $\gg(\ff)k[\ff]$
via the isomorphism $S\simeq k[\ff]$ is not contained in any minimal prime of $k[\ff]$. Therefore, by
torsionfreeness of $R$ over $k[\ff]$, the ideal $(\gg(\ff))R$ is not contained in any minimal prime of $R$.

(ii)
Let $S=k[\YY]/\mathfrak{b}$ be a common presentation of the algebras $k[\ff]$ and $k[\ff']$.
By assumption, there exist $F,F'\in R$ homogeneous such that $Ff'_j=F'f_j$, for $j=0,\ldots,m$ and, likewise,
$G=G(\yy),G'=G'(\yy)\subset S$ homogeneous such that $G(\yy)g'_{\ell}(\yy)=G'(\yy)g_{\ell}(\yy)$, for $\ell=0,\ldots,r$.
We have appended $\yy$ to the latter to emphasize that they are elements of $S=k[\YY]/\mathfrak{b}$.
Using the isomorphism $k[\YY]/\mathfrak{b}\simeq k[\ff']$, for instance, one has
\begin{equation}\label{eq_g}
G(\ff')g'_{\ell}(\ff')=G'(\ff')g_{\ell}(\ff').
\end{equation}
Substituting for $f'_j=(F'/F)f_j$, $j=0,\ldots,m$ in $g_{\ell}(\ff')$ and pulling out the obvious factor yields
$F^e\cdot g_{\ell}(\ff')=F'^e\cdot g_{\ell}(\ff)$, where $e$ stands for the common degrees of the $g_{\ell}$'s.
Multiplying by $F^e$ both members of (\ref{eq_g}) and
bringing into it the last equality gives
$$(F^e\,G(\ff'))g'_{\ell}(\ff')=G'(\ff')(F^e\, g_{\ell}(\ff'))= (G'(\ff')\,F'^e)g_{\ell}(\ff),$$
for $\ell=0,\ldots,r$.
Now, $G(\ff')$ and $G'(\ff')$ are regular (homogeneous) elements of $k[\ff] $ and $k[\ff']$, respectively as
they are both images of regular (homogeneous) elements of $S=k[\YY]/\mathfrak{b}$ under a graded $k$-isomorphism
of algebras. Since $R$ is torsionfree over $k[\ff']$ they are still regular (homogeneous) elements of $R$.
Therefore, $F^e\,G(\ff')$ and $G'(\ff')\,F'^e$ are regular homogeneous elements of $R$ and we are done.

{\sc Addendum}. As a slight control, if $d$ (respectively, $d'$) and $e$ (respectively, $e'$) are the respective degrees of
$f_j$ (respectively, $f'_j$) and of $g_{\ell}$ (respectively, $g'_{\ell}$) and, likewise
$D$ (respectively, $D'$) and $E$ (respectively, $E'$) are the respective degrees of
$F$ (respectively, $F'$) and of $G$ (respectively, $G'$), then
\begin{eqnarray*}
\deg ((G'(\ff')\,F'^e)+\deg (g_{\ell}(\ff))&=& d'E'+eD'+de= d'(E+e'-e) + e(D+d'-d)+de\\
& =& eD+d'E +d'e'
= \deg ((F^e\,G(\ff'))+\deg(g'_{\ell}(\ff')).\quad\quad\QED
\end{eqnarray*}

\subsubsection*{Image of a rational map}

 We have seen that all representatives of the same $m$-targeted rational map have a common
presentation $S=k[\YY]/\mathfrak{b}$, a standard graded reduced $k$-algebra.
Borrowing from the geometric terminology, we call Proj$(S)\subset \pp^m_k$ the {\sc image} of this map.

\subsubsection*{Restriction}

 One can introduce an additional notion, that of a restriction of an $m$-targeted rational map
$\mathfrak{F}:\ff\subset R$. Given a radical ideal $I\subset R$ such that $(\ff)$ is not contained in any minimal prime of $R/I$ on $R$
and such that $R/I$ is torsionfree over $k[\ff]/I\cap k[\ff]$, then we say that the map {\sc restricts} to Proj$(R/I)$.
Thus, being able to restrict a rational map requires quite a bit -- note that, if $I$ is a prime ideal, then
torsionfreeness of $k[\ff]/I\cap k[\ff]\subset R/I$ is trivially satisfied, so we just have to worry about the condition
on the relative minimal primes, which is unavoidable in the classical geometric setup.

\medskip

 There is one easy situation where restriction works as expected.

\begin{Lemma}\label{rational_vs_minimal_primes}
Let $\mathfrak{F}:\ff\subset R$ be an $m$-targeted rational map with source {\rm Proj}$(R)$.
\begin{enumerate}
\item[{\rm (1)}] If $P\subset R$ is a minimal prime of $R$ then
the map restricts to {\rm Proj}$(R/P)$ and  the image of the restriction is {\rm Proj}$(k[\ff]/P\cap k[\ff])$.
\item[{\rm (2)}] Conversely, if ${\wp}\subset k[\ff]$ is a minimal prime, then {\rm Proj}$(k[\ff]/{\wp})$
is the image of the restriction of $\ff$ to $R/P$ for some minimal prime $P\subset R$.
\end{enumerate}
\end{Lemma}
\demo
(1) This follows immediately from the properties of the rational map.

(2) It suffices to show that $\wp$ is contracted from some prime $\mathcal{P}\subset R_+$, as
in this case a minimal prime of $R$ contained in $\mathcal{P}$ will contract to $\wp$ as well.
For that, it suffices to show that $(\wp)R$ is not an irrelevant ideal in $R$.
But if that were the case, one would have $\wp=(\ff)k[\ff]$ which is impossible since $k[\ff]$
is a reduced ring of positive dimension (e.g., because $(\ff)$ has a regular element contained in
$k[\ff]$ -- see Remark~\ref{minors}).
\qed

\medskip

What the previous lemma does for us is to allow many arguments to be reduced to the case where $R$ is a integral domain.
Of course, this has to be used judiciously.

\begin{Remark}\rm (Geometric transcription)
Suppose that $k$ is algebraically closed.
By a version of the projective Nullstellensatz, (ii) implies that for every $P\in {\rm Min}(R)$ there exists some
homogeneous {\em submaximal} ideal (i.e., generated by $n-1$ linearly independent $1$-forms in $\xx$)
$\fn\subset R$ containing $P$ but not $(\ff)$.
For such an $\fn$, the contraction $\fn\cap k[\ff]\subset k[\ff]$ is an ideal generated by $m-1$  linearly independent
$1$-forms in $\ff$.
That is to say, if we write $k[\YY]/\mathfrak{b}\simeq k[\ff]$, then $\fn\cap k[\ff]$ is the residual image of
a homogeneous submaximal ideal in $k[\YY]$.
Geometrically, this implies that $f_0(\mathbf{a}),\ldots, f_m(\mathbf{a})$ makes sense as a point in $\pp^m$
for at least one point on any given irreducible component of the projective subset $V\subset \pp^n$ whose homogeneous
coordinate ring is $R$. Therefore, on every such component there is a nonempty open subset where $\ff$ affords a
well defined map.

In particular, the so obtained map $\mathfrak{F}$ is well defined on any open set $\mathcal{U}\subset V$ that intercepts non-emptily every irreducible component.
The image of $\mathfrak{F}$ is the closure $W\subset \pp^m$ of the set of images of points from any such $\mathcal{U}$.
This is the projective set of $\pp^m$ whose homogeneous coordinate ring is $k[\YY]/\mathfrak{b}$, i.e.,
after the obvious degree renormalization, $k[\ff]$ itself.
\end{Remark}

\subsubsection*{Birational maps and ``field'' extension criterion}

 The simplest instance of an $m$-targeted rational map is $\mathfrak{F}:\XX\subset R$,  represented by the residues in $R$ of the
variables $\{X_0,\ldots,X_n\}$ -- called {\sc the identity} rational map of Proj$(R)$.
The closest rational maps to the identity are the so-called birational maps.
We give a formal definition.

\begin{Definition}\label{def_of_birational}\rm
An $m$-targeted rational map $\ff\subset R=k[X_0,\ldots,X_n]/\mathfrak{a}$ of degree $d$ is said to be {\sc birational}
({\sc onto the image}) if there exists an $n$-targeted rational map $\mathfrak{G}:\gg\subset S=k[Y_0,\ldots,Y_m]/\mathfrak{b}\simeq k[\ff]$ such that
the composite $\gg(\ff)\subset R$ is equivalent to the identity on Proj$(R)$.
\end{Definition}
The rational map $\mathfrak{G}$ is called an {\sc inverse} to $\mathfrak{F}$.

The classical field theoretic criterion can be extended via the following formulation.

\begin{Proposition}\label{birational_in_terms_of_fields}
Let $\mathfrak{F}:\ff\subset R=k[X_0,\ldots,X_n]/\mathfrak{a}$ denote an $m$-targeted rational map of degree $d$.
The following conditions are equivalent:
\begin{enumerate}
\item[{\rm (a)}] $\mathfrak{F}:\ff\subset R$ is birational onto the image$\,${\rm ;}
\item[{\rm (b)}] For every minimal prime $P$ of $R$, the induced $m$-targeted rational map
$$\mathfrak{F}\restr P:\bar{\ff}\subset \bar{R}=R/P$$
 is birational onto its image and, moreover,
$R/P$ is the unique minimal component of $R$ mapping to the latter
by $\mathfrak{F}\,${\rm ;}
\item[{\rm (c)}] The induced homomorphism in degree zero $K(S)_{_0}\hookrightarrow K(R)_{_0}$  is an isomorphism$\,${\rm ;}
\item[{\rm (d)}] The induced homogeneous homomorphism $K(S)\hookrightarrow K(R^{(d)})$ is an isomorphism.
\end{enumerate}
\end{Proposition}
\demo The equivalence of (c) and (d) is fairly known -- see, e.g., \cite[Proof of Theorem 6.6]{ram2} or
\cite[Proof of Proposition 2.1]{SiVi}
for the argument in a more special setup.
Here we trade the non-homogeneous total ring of quotients for the homogeneous one, but the reasoning is the same, namely,
$K(S)=K(S)_{_0}(f)$ and $K(R^{(d)})=K(R^{(d)})_{_0}(f)=K(R)_{_0}(f)$ for some regular $1$-form $f\in S$ or, up to
identification $S=k[\ff]$,
a suitable $k$-linear combination of $f_{_0},\ldots,f_m$.
Since the maps in the statement are natural, the equivalence is clear.

Thus we are left with the mutual equivalence of (a), (b) and (c).

We note that, quite generally, if $R$ is any reduced Noetherian ring,
its ordinary total ring of quotients decomposes as the direct product
of the respective fields of fractions of the domains $R/P$, where $P$ runs through the minimal primes of $R$.
If, moreover, $R$ is standard graded over a field $k$ (possibly even more generally) then, by restriction, this decomposition
implies a similar decomposition of the {\em homogeneous} total ring of fractions
$K(R)=\prod _{P\in {\rm Min}(R)}K(R/P)$.
Clearly, this decomposition induces one in degree zero as well:
\begin{equation}\label{decomposition_degree0}
K(R)_{_0}=\prod _{P\in {\rm Min}(R)}K(R/P)_{_0},
\end{equation} where the factors
are now fields.

(a) $\Rightarrow$ (b)
This  piece follows from the first statement of Lemma~\ref{rational_vs_minimal_primes}.

(b) $\Rightarrow$ (c)
By the classical field criterion in the domain case one knows that he induced homomorphism in degree zero
$K(S/P\cap S)_{_0}\hookrightarrow K(R/P)_{_0}$  is an isomorphism for every minimal prime $P$ of $R$.
This is well-known, but for the sake of completion, we repeat the argument.
In fact, it is now easy to work directly with coordinates.
Updating the notation, assume that $R=k[\xx]$ is a domain.
Then one has
$$K(S)_{_0}=k(\yy)_{_0}=k\left(\frac{y_1}{y_0},\ldots,\frac{y_m}{y_0}\right)\simeq k(\ff)_{_0}=k\left(\frac{f_1}{f_0},\ldots,\frac{f_m}{f_0}\right)
\subset k(\xx)_{_0}=k\left(\frac{x_1}{x_0},\ldots,\frac{x_n}{x_0}\right),$$
where we assume $y_0\neq 0, x_0\neq 0$ without lost of generality.
By assumption, there are nonzero (homogeneous) $H,H'\in R$ such that $H\cdot g_{\ell}(\ff)=H'\cdot x_{\ell}$, for $\ell=0,\ldots,n$.
Then
$$\frac{x_{\ell}}{x_0}=\frac{g_{\ell}(\ff)}{g_0(\ff)}=\frac{f_0^e\left(g_{\ell}(f_1/f_0,\ldots,f_m/f_0)\right)}{f_0^e\left(g_{_0}(f_1/f_0,\ldots,f_m/f_0)\right)}
=\frac{g_{\ell}(f_1/f_0,\ldots,f_m/f_0)}{g_{_0}(f_1/f_0,\ldots,f_m/f_0)},\quad \ell=0,\ldots,n,$$
thus showing the required equality $k(f_1/f_0,\ldots,f_m/f_0)=k(x_1/x_0,\ldots,x_n/x_0)$.

Next, note that the map $K(S/P\cap S)_{_0}\hookrightarrow K(R/P)_{_0}$ is naturally induced by the map $K(S)_{_0}\hookrightarrow K(R)_{_0}$.

Then, up to identification we only need to show that if $K(S/S\cap P)_{_0}=K(R/P)_{_0}$ for every $P\in {\rm Min}(R)$, then $K(S)_{_0}=K(R)_{_0}$.
To see this, note that the hypothesis implies that, for every $\wp\in {\rm Min}(S)$, exactly one prime $P\in {\rm Min}(R)$
contracts to $\wp$ -- i.e., the map Min$(R)\rar {\rm Min}(S)$ is surjective (hence, bijective by the second
statement of Lemma~\ref{rational_vs_minimal_primes}).
This is due to the fact that in the injection $K(S)_{_0}\hookrightarrow K(R)_{_0}$ the decomposition (\ref{decomposition_degree0})
forces the restriction to an individual factor $K(S/\wp)$ to map diagonally to the product of the factors of $K(R/P)_{_0}$, one for each
$P\in {\rm Min}(R)$ contracting to $\wp$.
By the same token, this implies the sought global isomorphism $K(S)_{_0}\simeq K(R)_{_0}$.

(c) $\Rightarrow$ (a)  (Assuming $k$ infinite)
Fix a regular homogenous element $x\in (\xx)=R_+$ of the form $x_0+\sum_{\ell =1}^n\lambda_{\ell}x_{\ell}$, for suitable $\lambda_{\ell}\in k$.
For $\ell=1,\ldots,n$, let $g_{\ell}(\yy)/g(\yy)\in K(S)_{_0}$
denote the (uniquely defined) inverse image of $x_{\ell}/x$ under the isomorphism $K(S)_{_0}\simeq K(R)_{_0}$, for suitable
homogeneous elements $g(\yy),g_{\ell}(\yy)\in (\yy)=S_+$, with $g(\yy)$ regular.
Since the isomorphism $K(S)_{_0}\simeq K(R)_{_0}$ is induced by $y\mapsto f_j$, one has $g_{\ell}(\ff)/g(\ff)=x_{\ell}/x$, for $\ell=1,\ldots, n$.
Thus, $x\cdot g_{\ell}(\ff)=g(\ff)\cdot x_{\ell}$, for $\ell=1,\ldots, n$, with $x, g(\ff)$ regular homogeneous elements in $R$.
Now set $g_0(\yy):=g(\yy)-\sum_{\ell =1}^n\lambda_{\ell}g_{\ell}(\yy)$.
Then
\begin{eqnarray*}g(\ff)x_0&=&g(\ff)\left(x-\sum_{\ell =1}^n\lambda_{\ell}x_{\ell}\right)=xg(\ff)-\sum_{\ell =1}^n\lambda_{\ell}(g(\ff)x_{\ell})=
x g(\ff)-\sum_{\ell =1}^n\lambda_{\ell}(x g_{\ell}(\ff))\\
&=&x\left(g(\ff)-\sum_{\ell= 1}^n\lambda_{\ell}g_{\ell}(\ff)\right)=x\cdot g_0(\ff).
\end{eqnarray*}
Finally, note that $g_0(\yy),\ldots, g_n(\yy)\subset S=k[\yy]$ is a rational $(n+1)$-datum by Lemma~\ref{stability}, hence
defines an $n$-targeted rational map.
Therefore, the rational $(n+1)$-data $g_0(\ff),\ldots, g_n(\ff)$ and $x_0,\ldots,x_n$ are equivalent, as was to be shown.
\qed

\begin{Corollary}\label{reverse_composition}
If an $m$-targeted rational map $\mathfrak{F}:\ff\subset R=k[X_0,\ldots,X_n]/\mathfrak{a}$ is birational
{\rm (}onto its image{\rm )} with $n$-targeted inverse rational map $\mathfrak{G}:\gg\subset S=k[Y_0,\ldots,Y_m]/\mathfrak{b}\simeq k[\ff]$,
then the composite $\ff(\gg)\subset S$ is equivalent to the identity on Proj$(S)$.
Consequently, the inverse is uniquely defined.
\end{Corollary}
\demo The first statement follows from the equivalence of (a) and (c) in the previous theorem.
The second statement follows from the first in the usual way, namely writing $\equiv$ for equivalence, if $\gg'$
is another inverse, we have
$$\XX\equiv\gg(\ff) \Rightarrow  \gg'\equiv\XX (\gg')\equiv \gg(\ff)(\gg')\equiv\gg(\ff(\gg'))\equiv
\gg(\YY)\equiv \gg.$$
\qed

The corollary means, as it should, that a birational map on Proj$(R)$ onto its image Proj$(S)$ has a
unique ``two-sided'' inverse and the inverse is a birational map on Proj$(S)$ onto Proj$(R)$.

\subsection{The JD-rank and its bearing to birationality}

\subsubsection*{Definition and background}

Let us recap the present algebraic tool for birationality that
goes in parallel to the geometric notion of the linear fibers of a
rational map.
This concept has evolved continuously from previous notions, of
which the cradle is the Jacobian dual fibration of \cite{dual} (see
the subsequent \cite{cremona} and \cite{bir2003}). In order to avoid
tedious repetition, we take a slightly diversion by looking into
ideals (respectively, codimension) instead of matrices
(respectively, rank).

As we see, the base ideal of a rational map will play an essential role not entirely
foreseen in the previous subsection.

Thus, let $\mathfrak{F}:\ff\subset R=k[\XX]/\mathfrak{a}$ be an $m$-targeted rational map with source
Proj$(R)$, where $R$ is as in the previous subsection, with $\dim R\geq 1$.

Since $\ff$ is generated in same degree, the Rees algebra ${\cal R}_R((\ff))$ is a bigraded
standard $k$-algebra. In particular, there is a presentation $R[\YY]=R[Y_0,\ldots,Y_m]\surjects {\cal R}_R((\ff))$
mapping $Y_j\mapsto f_j$.
Let  ${\cal J}=\ker(R[\YY]\surjects {\cal R}_R((\ff))$ denote the corresponding kernel, often
called a presentation ideal of ${\cal R}_R(\ff)$.
Just as the image of a rational map, $\mathcal{J}$ depends only on the rational map
and not on any of its representatives.
Note that it is a bihomogeneous ideal
in the natural standard bigrading of $R[\YY]$ induced by the one of the
bivariate polynomial ring $k[\XX,\YY]$. Thus, one has
$${\cal J}=\bigoplus_{(p,q)\in \mathbb{N}^2} {\cal J}_{(p,q)},$$
where ${\cal J}_{(p,q)}$ denotes the $k$-vector space of forms of bidegree $(p,q)$.

The first relevant piece for our purpose is ${\cal J}_{0,*}$, which is
spanned by the forms of bidegree  $(0,q)$, for all $q\geq 1$.
Since these forms have coefficients only in $k$, one can view them as elements of $k[\YY]$.
As such they generate the ideal of all polynomial relations (over $k$) of $\ff$.
Thus,
$$k[\YY]/({\cal J}_{0,*})k[\YY]\simeq k[\ff],$$
as standard graded $k$-algebras,
after a degree renormalization in $k[\ff]$.
Thus, the ideal $({\cal J}_{0,*})k[\YY]$ recovers the ideal $\mathfrak{b}$ of the image brought up in the previous section.

We write $S=k[\YY]/({\cal J}_{0,*})k[\YY]\simeq k[\ff]$ throughout, thinking interchangeably of $S$ as one or the other.

\medskip

For birationality, the following bihomogeneous piece is as important:
$${\cal J}_{1,*}:=\bigoplus_{r\in\NN} {\cal J}_{1,q}$$
with ${\cal J}_{1,q}$ denoting the bigraded piece of ${\cal J}$ spanned by the forms of bidegree
 $(1,q)$ for all $q\geq 0$. Note that ${\cal J}_{1,*}=\{0\}$
 may very well be the case in general.

Now, a form of bidegree $(1,*)$ can be written as $\sum_{i=0}^n Q_i(\YY)\,x_i$, for suitable homogeneous $Q_i(\YY)\in k[\YY]\subset R[\YY]$
of the same degree.
Since $\YY$ are indeterminates over $R$, two such representations of the same form imply a syzygy of $\{x_0,\ldots,x_n\}$
with coefficients in $k$. Thus the representation is unique  up to $k$-linear dependency
of $\{x_0,\ldots,x_n\}$,  i.e., up to elements of ${\mathfrak a}_1$.

In particular, if the {\sc embedding dimension} of $R$ -- i.e., the $k$-vector space dimension
$\edim (R):=\dim_k(R_1)=n+1-\dim_k {\mathfrak a}_1$ -- is $n+1$ then every such form has a unique expression.

Next, one can pick a minimal set of generators of the ideal $({\cal J}_{1,*})$ consisting of a finite number
of forms of bidegree $(1,q)$, for various $q$'s.
Let $\{P_1,\ldots,P_s\}\subset k[\XX,\YY]$ denote liftings of these biforms and let $\{\ell_1,\ldots,\ell_r\}\subset k[\XX]$ be  a $k$-vector space
basis of ${\mathfrak a}_1$.
Consider the Jacobian matrix of the polynomials $\{\ell_1,\ldots,\ell_r, P_1,\ldots,P_s\}$ with respect to $\XX$, a matrix with entries in $k[\YY]$.
Write $\psi$ for the corresponding matrix over $S=k[\YY]/{\mathfrak b}$, here called the weak Jacobian matrix associated to
the given set of generators of $({\cal J}_{1,*})$.
Note that $\psi$ has a distinguished $(n+1-\edim(R))\times (n+1)$ submatrix $\rho$ with entries in $k[\YY]_0=k$,
which is non-singular. Its number of rows $r=n+1-\edim(R)$ can be dubbed the {\sc degeneration index} dgi$(\mathfrak{a})$ of the
presentation $R\simeq k[\XX]/\mathfrak{a}$.
Clearly, this index is a lower bound for the rank of $\psi$ over $S$ whenever this rank is well-defined.

On the other hand, by definition one clearly has $\XX\cdot \psi^t=0$ over ${\cal R}_R(\ff)$,
in particular the vector $(\XX)^t$ belongs to the null space of $\psi$ over ${\cal R}_R(\ff)$.
If $R$ is a domain then the rank of $\psi$ over ${\cal R}_R(\ff)$ is at most $n$, hence for even more reason
also over the residue ring ${\cal R}_R((\ff))\surjects S$.

A consequence of Proposition~\ref{rank_and_dimension} below will show that this upper bound holds true even
if $R$ is not integral, provided $\psi$ still has a rank.

Note that a weak Jacobian matrix $\psi$ is not uniquely defined due to the lack of uniqueness in the expression of
an individual form and to the choice of bihomogeneous generators.
Besides, it  can be fairly ``slack'', in the sense that some of its
columns may be zero -- actually $\psi$ itself can turn out to be the zero matrix for that matter!

On the bright side however, one has:

\begin{Lemma}\label{rank_is_stable}
If the weak Jacobian matrix associated  to a set of bihomogeneous minimal generators of
$({\cal J}_{1,*})$ has rank over $S$ then the weak Jacobian matrix associated to any other
set of bihomogeneous minimal generators of
$({\cal J}_{1,*})$ has  rank over $S$ and the two ranks coincide.
\end{Lemma}
\demo
One shows, namely, that two such matrices have the same column submodule over $S$.

Let $\{P_1,\ldots, P_s\}, \{Q_1,\ldots, Q_s\}\subset k[\XX,\YY]$ denote liftings of two sets of minimal
biforms of various bidegrees  $(1,q)$ , for suitable $q\geq 1$, generating $({\cal J}_{1,*})$.
Note that though the sets may be different, the set of corresponding bidegrees is uniquely defined up
to reordering of generators.
Say, $\deg (P_t)=\deg (Q_t)=(1,d_t)$, with $1\leq t\leq s$ and $1\leq d_1\leq d_2\leq\cdots\leq d_s$.
Therefore, for $t=1,\ldots,s$ we may write
$$Q_t=\sum_{j=1}^s \lambda_{tj} p_{tj}(\YY)P_j + \sum_{|\alpha_t|=d_t} a_{\alpha_t}(\XX)\, \YY^{\alpha_t},$$
for suitable $\lambda_{tj}\in k$, $p_{tj}(\YY)\in k[\YY]$ of degree $d_t-\deg_{\YY}(P_j)$ and
$a_{\alpha_t}(\XX)\in{\mathfrak a}_1$.
Setting $a_{\alpha_t}(\XX)=\sum_{k=1}^r \mu_{tk}\ell_k(\XX)$, we obtain
\begin{equation*}
\frac{\partial Q_t}{\partial X_i}=\sum_{j=1}^s \, \lambda_{tj} p_{tj}(\YY)\frac{\partial P_j}{\partial X_i} +
\sum_{|\alpha_t|=d_t} \left(\sum_{k=1}^r \mu_{tk}\,\YY^{\alpha_t}\frac{\partial \ell_k(\XX)}{\partial X_i} \right)
\end{equation*}
This shows that the $t$th column of the Jacobian matrix of $\{\ell_1,\ldots,\ell_r, Q_1,\ldots,Q_s\}$
is contained in the column $k[\YY]$-submodule of the Jacobian matrix of $\{\ell_1,\ldots,\ell_r, P_1,\ldots,P_s\}$.
By symmetry, the two column $k[\YY]$-submodules coincide, hence the $S$-submodule Im$(\psi)$ through any one of
the two sets of generators is the same, as was to be shown.
\qed

\medskip

Because of the previous lemma, we say that any weak Jacobian matrix associated to a set of minimal generators of $({\cal J}_{1,*})$
is a {\sc Jacobian dual matrix} of $\ff$ (or of the corresponding $m$-targeted rational map $\mathfrak{F})$.
The lemma unveils a numerical character of such a construction.

\begin{Definition}\label{jacobian_rank}\rm
If some (hence, every) Jacobian dual matrix of an $m$-targeted rational map $\mathfrak{F}:\ff\subset R\simeq
k[\XX]/\mathfrak{a}$ has a rank over $S$, we will say that $\mathfrak{F}$  has a {\sc Jacobian dual rank}
and write ${\rm jd}\rk(\mathfrak{F})$  to denote it.
In this situation, the {\sc non-degenerate Jacobian dual rank} of $\mathfrak{F}$ is ${\rm jd}\rk_+(\mathfrak{F}):
={\rm jd}\rk(\mathfrak{F})-{\rm dgi}(\mathfrak{a})\geq 0$.
\end{Definition}
Note that the non-degenerate Jacobian dual rank is well-defined, provided the Jacobian dual rank is.

The next result shows that the rank just introduced is sensitive to the dimension
difference between source and target.

\begin{Proposition}\label{rank_and_dimension}
Let $\mathfrak{F}:\ff\subset R$ stand for an $m$-targeted rational map.
Set $S=k[\YY]/\mathfrak{b}\simeq k[\ff]\subset R$.
If $\mathfrak{F}$ has a Jacobian dual rank
then
$$\dim R - \dim S \leq n - {\rm jd}\rk(\mathfrak{F})=\edim(R)-1-{\rm jd}\rk_+(\mathfrak{F}).$$
\end{Proposition}
\demo Up to a non-singular $k$-linear map on $R$, one may assume that the forms in $\ff$ involve effectively at most
the first $r+1=\edim(R)$ elements $x_0, \ldots, x_r$.
More precisely, let
$$A=\left(
\begin{array}{c@{\quad\vrule\quad}c}
\mathbf{1}_{r+1} & \mathbf{a}=(a_{h,i}) \\
\multispan2\hrulefill\\
  & \mathbf{1}_{n-r}
\end{array}
\right)$$
be an invertible $k$-matrix such that $(X_0,\ldots,X_n)\cdot A =(X_0,\ldots,X_r,\ell_{r+1},\ldots, \ell_n)$,
with $\ell_h=X_h-\sum_{i=0}^r a_{h,i}X_i, r+1\leq h\leq n$
and $\{\ell_{r+1},\ldots, \ell_n\}$ forming
a basis of the $k$-vector space $\mathfrak{a}_1$.
This defines a  homogeneous $k$-algebra isomorphism
$$R=k[X_0,\ldots,X_n]/\mathfrak{a}\simeq
\tilde{R}=k[X_0,\ldots,X_r]/\tilde{\mathfrak{a}}$$
which in turn restricts to homogeneous $k$-algebra isomorphism of the respective subalgebras
$k[\ff]$ and $k[\tilde{\ff}]$, where $f_j\mapsto \tilde{f_j}$ under this map --
however, note that the last isomorphism is not induced by equivalence of representatives, one is rather changing
the rational map by ``source and target'' transformations.

Thus,  this operation does not affect either $\dim R$ or $\dim S$ and, likewise, will not affect  the numbers in
the rightmost side of the stated inequality.

Updating the notation, consider the $S$-module $E:={\rm coker}_S({\psi}^t)$ defined as the cokernel of the dual map associated to ${\psi}$.
Since ${\psi}^t$ has a rank over $S$, so does $E$. Then
$$\dim \mathcal{R}_S(E)=\dim S+\rk_S(E)=\dim S+\edim(R)-\rk_S(\psi),$$
where $\mathcal{R}_S(E)$ stands for the Rees algebra of $E$ in the sense of \cite{ram1}).
On the other hand, a piece of the assumption on the base ideal $(\ff)$ implies that
$\dim \mathcal{R}_R((\ff))=\dim R+1$ -- for which it suffices to know that $(\ff)$ is not contained in at least one minimal prime of
$R$ of maximal dimension.

Then the asserted inequality follows from
the inequality
$$\dim \mathcal{R}_R((\ff))\leq \dim \mathcal{R}_S(E),$$
which in turn stems from the existence of an $S$-algebra surjection $\mathcal{R}_S(E)\surjects \mathcal{R}_R((\ff))$,
with $\mathcal{R}_R((\ff))$ viewed as an $S$-algebra via
$$S=k[\YY]/\mathfrak{b}\simeq k[\YY]/({\cal J}_{0,*})\injects R[\YY]/\mathcal{J},$$
where $\mathcal{J}$ is the presentation ideal of $(\ff)$ based on these generators.

To see this, one first shows that there is an $S$-surjection $\mathcal{S}_S(E)\surjects \mathcal{R}_R((\ff))$ with source the
symmetric algebra of $E$. Indeed,  one has
\begin{eqnarray*}
\mathcal{S}_S(E)&\simeq &S[{\XX}]/I_1({\XX}\cdot{\psi}^t)
=k[\YY][{\XX}]/(\mathfrak{b},I_1({\XX}\cdot{\psi}^t))\surjects {R}[\YY](\mathfrak{b},I_1({\XX}\cdot{\psi}^t))\\
&\surjects &   R[\YY]/\mathcal{J}.
\end{eqnarray*}

According to \cite{ram1}, $\mathcal{R}_S(E)\simeq \mathcal{S}_S(E)/(0:_{\mathcal{S}_S(E)} s),$
the symmetric algebra modulo torsion, where $s$ is a suitable regular homogeneous element of $S$.
Therefore, it suffices to show that $0:_{\mathcal{S}_S(E)} s$  maps to zero under the $S$-surjection
$\mathcal{S}_S(E)\surjects \mathcal{R}_R((\ff))$.
Thus, pick an element $G\in 0:_{\mathcal{S}_S(E)} s$ that maps to a nonzero element of $\mathcal{R}_R((\ff))$.
Then $s\in S\subset \mathcal{R}_R((\ff))$ under the above injection annihilates a nonzero element of $\mathcal{R}_R((\ff))$.
Since $\mathcal{R}_R((\ff))$ is $R$-torsionfree and $S$-torsionfree as a piece of the assumption that the inclusion $S\subset R$ is induced
by a rational map, then $\mathcal{R}_R((\ff))$ is $S$-torsionfree.
Therefore, we get a contradiction.
\qed

\begin{Corollary}\label{maximal_rank_of_psi}
With the notation of the previous theorem, ${\rm jd}\rk_+(\mathfrak{F})\leq \edim(R)-1$.
Moreover, if equality holds then $\dim S=\dim R$.
\end{Corollary}

Thus, if the rank exists, one has the bounds:
\begin{eqnarray*}
&&{\rm dgi}(\mathfrak{a})\leq {\rm jd}\rk(\mathfrak{F})\leq n\\
&& 0\leq {\rm jd}\rk_+(\mathfrak{F})\leq \edim(R)-1.
\end{eqnarray*}

\begin{Remark}\rm We note that, in the classical integral case in characteristic zero, the corollary says:
${\rm jd}\rk_+(\mathfrak{F})= \edim(R)-1$ implies that the general fiber of the map
is a finite set of points.
Thus, the preliminary estimate in Proposition~\ref{rank_and_dimension} puts in evidence
that the Jacobian dual rank
 plays a sort of complementary role to a transcendence degree in
the integral case.
The subsequent results will actually show that it takes over a complementary role to the
degree of the rational map in the classical setup.
\end{Remark}

\subsection{Statement and proof of the criterion}

We now state the main theorem of this section, which gives a criterion of birationality.
After having worked all this way a purely algebraic setup, it is natural to first convey an algebraic formulation
of the criterion.

We keep the essential notation of the previous subsections.

\begin{Theorem}\label{main_criterion_algebraic}
Let $\mathfrak{F}:\ff\subset R=k[\XX]/\mathfrak{a}=k[X_0,\ldots,X_n]/\mathfrak{a}$ stand for an $m$-targeted rational map and,
for every minimal prime $P$ of $R$, let $\mathfrak{F}\restr P$ denote the induced $m$-targeted rational map with source
$R/P$.
The following conditions are equivalent:
\begin{itemize}
\item[{\rm (a)}] $\mathfrak{F}$ is birational onto the image.
\item[{\rm (b)}] $\mathfrak{F}$ has non-degenerate Jacobian dual rank satisfying
$${\rm jd}\rk_+(\mathfrak{F})= \edim(R)-1.$$

\item[{\rm (c)}] For every minimal
prime $P$ of $R$, $\mathfrak{F}\restr P$ has non-degenerate Jacobian dual rank satisfying
$${\rm jd}\rk_+(\mathfrak{F}\restr P)= \edim(R/P)-1.$$

\end{itemize}
{\sc Supplement.} Suppose, moreover, that the above equivalent conditions holds.
Let $\psi$ denote a weak Jacobian dual matrix of $\mathfrak{F}$ over $S=k[\YY]/\mathfrak{b}\simeq k[\ff]$.
\begin{enumerate}
\item[{\rm (i)}] A rational $(n+1)$-datum in
$S$ representing the inverse map is given by the coordinates of any homogeneous vector of positive degree
in the {\rm (}rank one{\rm )} null space of $\psi$ over $S$ for which these coordinates generate an ideal containing
a regular element.

\item[{\rm (ii)}] If further $R$ is a domain, a rational $(n+1)$-datum as in {\rm (i)}
can be taken to be the set of the {\rm (}ordered, signed{\rm )} $(\edim(R)-1)$-minors
of an arbitrary $(\edim(R)-1)\times \edim(R)$ submatrix of $\psi$ of rank $\edim(R)-1$.
\end{enumerate}
\end{Theorem}
\demo
As in the proof of Proposition~\ref{rank_and_dimension}, up to a non-singular $k$-linear map on $R$, one may assume that the forms
in $\ff$ involve only the first $r+1=\edim(R)$ elements $x_0, \ldots, x_r$.
In addition, the resulting isomorphisms at the respective levels of $R$ and $S$ induce homogeneous
$k$-algebra isomorphisms $K(R)\simeq K(\tilde{R})$ and $K(S)\simeq K(\tilde{S})$, where
as before $S=k[\YY]/\mathfrak{b}\simeq k[\ff]\simeq k[\tilde{\ff}]$.

Finally, passing to degree zero and using the field criterion of Proposition~\ref{birational_in_terms_of_fields},
it follows that $\mathfrak{F}:\ff\subset R$ is birational with source Proj$(R)$ if and only if $\tilde{\ff}\subset \tilde{R}$
gives a birational map with source Proj$(\tilde{R})$.
This shows that (a) is stable under this change.
According to Proposition~\ref{rank_and_dimension} and the discussion of its preliminaries, the
non-degenerate Jacobian dual rank in both (a) and (b) is also stable under this process.

Thus, we may restart now assuming that $n+1=\edim(R)$ -- in other words, $R=k[X_0,\ldots, X_n]/\mathfrak{a}$, with
$\mathfrak{a}\subset (R_+)^2$ -- and ${\rm jd}\rk(\mathfrak{F})={\rm jd}\rk_+(\mathfrak{F})$.

\medskip

(a) $\Rightarrow$ (b)
By assumption, there is a rational $(n+1)$-datum $\gg(\yy)=\{g_0(\yy),\ldots, g_n(\yy)\}\subset S$ such that $\gg(\ff)$ is equivalent
to the identity on $R$.
Now, there is an isomorphism of $k$-algebras $\Phi:\mathcal{R}_R ((\ff))\simeq \mathcal{R}_S((\gg))$ induced by the identity map of
$k[\XX,\YY]/(\mathfrak{a},\mathfrak{b})$. The existence of this isomorphism is the same as in \cite[Theorem 2.1, (i) $\Rightarrow$ (ii)]{bir2003}, observing that the hypothesis of $R$ (hence $S$) being a domain is superfluous because the two maps constructed in [loc.cit.] were inverse of each other by construction.
We repeat the basic idea for the sake of completeness and further clarity while adapting the argument to the present context.
One has isomorphisms
$$\frac{k[\xx,\yy]/(\mathfrak{a},\mathfrak{b})}{{\cal J_{\ff}}}\simeq \frac{k[\xx]}{\mathfrak{a}}\,[\ff t]=\mathcal{R}_R ((\ff)),\;\;
\mbox{\rm mapping}\;\;x_i\mapsto x_i (\bmod \mathfrak{a}),\; y_j\mapsto f_j t$$
and, similarly,
$$\frac{k[\xx,\yy]/(\mathfrak{a},\mathfrak{b})}{{\cal J_{\gg}}}\simeq \frac{k[\yy]}{\mathfrak{b}}\,[\gg u]=\mathcal{R}_S((\gg)),\;\;
\mbox{\rm mapping}\;\; y_i\mapsto y_i (\bmod \mathfrak{b}),\; x_i\mapsto
g_i u.$$
Then $\Phi$ drives ${\cal J_{\ff}}$ inside ${\cal J_{\gg}}$ by showing that if $F(\xx,\yy)\in k[\xx,\yy]$ is
bihomogeneous such that  $F(\xx,\ff t)\equiv 0 \pmod {\mathfrak{a}}k[\xx,t]$ then
$F(\gg u,\yy)\equiv 0  \pmod {\mathfrak{b}}k[\yy,u]$.
To prove this one uses that the composite $\ff(\gg)$ is equivalent to the identity
of Proj$(S)$. The reverse map works the same way by the obvious exchange.

Restricting to the bihomogeneous pieces, one has in particular $\Phi((\mathcal{J}_{\mathbf{f}})_{(1,s)})=(\mathcal{J}_{\mathbf{g}})_{(s,1)}$,
for every $s\geq 1$.
Note that for varying $s\geq 1$, $(\mathcal{I}_{\mathbf{g}})_{(s,1)}$ generate
the presentation ideal of the symmetric algebra of $\mathbf{g}$ over $S$.
Up to a permutation of the variables $\xx$, one has, say, $(\mathcal{I}_{\mathbf{g}})_{(s,1)}=I_1((x_0\cdots x_{\mu-1})\cdot \mathcal{M})$,
where $\mathcal{M}$ denotes the matrix of a
 minimal graded presentation of $\mathbf{g}$ over $S$ and
$\mu=\mu(\mathbf{g})$ is the minimal number of generators of the ideal $(\gg)$.
Further, if $\psi$ denotes a weak Jacobian dual matrix of $\mathbf{f}$ over $S$ then
$\Phi((\mathcal{I}_{\mathbf{f}})_{(1,s)})=I_1((x_0\cdots x_n)\cdot \psi^t)$ by the definition
of $\psi$.
The ensued equality $I_1((x_0\cdots x_n)\cdot \psi^t)=I_1((x_0\cdots x_{\mu-1})\cdot \mathcal{M})$
shows that the last $n-\mu$ rows of $\psi^t$ are null.
Hence if $\psi'$ is the matrix obtained by removing the zero columns of $\psi$ then
any column of $\psi'^t$ lives in Im$(\mathcal{M})$ and vice versa.
Therefore Im$(\psi')=$ Im$(\mathcal{M})$.

Now $\mathcal{M}$ has a  rank over $S$. This is because $\grade_S(\mathbf{g})\geq 1$ by by definition of a rational datum,
hence a presentation map of $(\gg)$ over $S$ has a rank as well.
Clearly then $\psi'$ has a rank  and $\rk(\psi)=\rk(\psi')=\rk(\mathcal{M})=\mu-1=\edim(R)-1$,
the last equality following from the non-degeneracy of $R$.

\medskip

(b) $\Rightarrow$ (c)
Note that $n=\edim(R)-1$ as we are assuming that $R$ is non-degenerate.
By the assumption and Lemma~\ref{rank_and_dimension}, the ideal of minors $I_n(\psi)$
of a weak Jacobian dual matrix $\psi$ has a regular element on $S$.
Let $P$ be minimal prime of $R$ and $Q=P\cap S$, using the isomorphism $S\simeq k[\ff]$.
Clearly then $I_n(\psi\restr Q)\neq \{0\}$, where $\psi\restr Q$ is the matrix
on $S/Q$ obtained by taking the entries of $\psi$ modulo $Q$.
It follows that $\psi\restr Q$ has rank $n$ over $S/Q$.

Now consider a weak Jacobian dual matrix $\tilde{\psi}$ of $\mathfrak{F}\restr P$, as obtained from the Rees
algebra $\mathcal{R}_{R/P}((\ff)+P/P)$. By definition, Im$(\psi\restr Q$) can be taken to
be a submodule of Im$(\tilde{\psi})$.
Therefore, the rank of $\tilde{\psi}$ is at least $n$, hence it obtains
$$\edim(R/P)-1 - {\rm jd}\rk_+(\mathfrak{F}\restr P)=n- {\rm jd}\rk(\mathfrak{F}\restr P)\leq 0,$$
from which ensues ${\rm jd}\rk_+(\mathfrak{F}\restr P)\geq \edim(R/P)-1$
and then equality follows by Proposition~\ref{rank_and_dimension}, as required.

\medskip

(c) $\Rightarrow$  (a)
The proof will simultaneously argue for item (ii).

The first step is to use Proposition~\ref{birational_in_terms_of_fields}, by which the condition in (a) is equivalent
to the condition that every restriction $\mathfrak{F}\restr P$ of $\mathfrak{F}$ is birational onto its image
and, moreover, $R/P$ is the unique minimal component of $R$ mapping to the latter
by $\mathfrak{F}\,$.

To show that one can reduce to the case where $R$ is a domain, we argue as follows.
Given a minimal prime $P$ of $R$, write
$R/P\simeq (k[\xx]/\mathfrak{a})/(\mathfrak{P}/\mathfrak{a})\simeq k[\xx]/\mathfrak{P}$
and $S/P\cap S\simeq k[\yy]/\mathfrak{Q}$, where $\mathfrak{Q}$ is the inverse image of
$\mathfrak{P}$ in $k[\yy]$.
Further, set $\mathfrak{F}_{\mathfrak{P}}=\mathfrak{F}\restr P$ to emphasize the original polynomial
ambient.
The assumption then reads as
$${\rm jd}\rk_+(\mathfrak{F}_{\mathfrak{P}})= \edim(R/P)-1=n-{\rm dgi}(\mathfrak{P}),$$
where dgi$(\mathfrak{P})$ is the degeneration index of $\mathfrak{P}$ as discussed in the beginning of the subsection.
Equivalently, ${\rm jd}\rk(\mathfrak{F}_{\mathfrak{P}})=  n$.
Thus, if we show that this implies that $\mathfrak{F}_{\mathfrak{P}}$ is birational onto its image
we will be done.

Therefore, we may assume at the outset that $R$ is a domain and a weak Jacobian dual matrix
$\psi$ has rank $n$.
But then one can select an $n\times (n+1)$ submatrix $\sigma$ of $\psi$ of rank $n$.
Let ${\bf \Delta}:=\{\Delta_0,\ldots ,\Delta_n\}$ stand for the ordered signed $n$-minors of $\sigma$.
This set is a rational $(n+1)$-datum.
It remains to prove that $\mathbf{\Delta}(\ff)$ is equivalent to the identity on Proj$(R)$.

Apply Proposition~\ref{koszul_hilbert} with $A:=k[\yy]$ and with $\rho$ standing
for the transpose of the selected $n\times (n+1)$ submatrix $\sigma$ of $\psi$ of rank $n$ to get that
the Koszul forms $\Delta_i\,x_j-\Delta_j\,x_i$, one for each pair $i,j$, belong to the
$k[\yy]$-submodule of $A^{n+1}=\sum_{i=0}^{n}\,A\, x_i$ generated by the forms induced by the rows of $\sigma$
-- i.e, obtained as inner product of these by the transpose $(x_0,\ldots,x_n)^t$.
But the latter, by definition, belong to the defining ideal of the Rees algebra
of $(\ff)$ over $R[\YY]$, hence vanish on the tuple whose coordinates are the constituents of $\ff$.

Therefore
\begin{equation}\label{koszul}
\Delta_i(\ff)\,x_j=\Delta_j(\ff)\,x_i, \quad \mbox{\rm for each pair} \;i,j.
\end{equation}

Recall that some $\Delta_i(\ff)\neq 0$; say, this is the case for $i=n$.
On the other hand,
$x_n\neq 0$ since otherwise (\ref{koszul}) implies that $x_i=0$ for every $i$, which is nonsense.
In particular, one has the equalities $\Delta_n(\ff)\, x_i=x_n\, \Delta_i(\ff)$, with $\Delta_n(\ff), x_n$
nonzero.
Therefore $\mathbf{\Delta}(\ff)$ is equivalent to the identity on Proj$(R)$.

\medskip

It remains to prove (i).

Let $\hh^t\in\ker(\psi^t)$ be given as in the hypothesis of (i).
Let $\gg$ be a representative of the inverse map such that the vector $\gg^t\in\ker(\psi^t)$
as in the proof and notation of the implication (a) $\Rightarrow$ (b).
The argument there showed that $\ker(\psi^t)$ has rank $1$.
Since the ideal $(\hh)$ has a regular element by assumption, Lemma~\ref{minors} implies
that $\hh$ and $\gg$ are equivalent $(n+1)$-sets.
Then, by Lemma~\ref{stability}, $\hh$ is rational datum, hence defines the inverse map as well.
\qed

\begin{Remark}\rm The argument in the proof of (ii) above goes through regardless of the integrality
assumption, provided a weak jacobian dual matrix of $\ff$ admits a  submatrix of size
$(\edim(R)-1)\times \edim(R)$ and of rank $\edim(R)-1$.
\end{Remark}

We now transcribe the above result into a geometric version.
We assume that the ground field $k$ is infinite and that varieties
have points in an algebraic closure of $k$.

\begin{Theorem}\label{main_criterion_geometric}
Let $V\subset \pp^n\, (n\geq 1)$ denote a reduced, but not necessarily irreducible,
projective subvariety.
Let $\mathfrak{F}:V\dasharrow \pp_k^m$ denote a rational map.

The following conditions are equivalent:
\begin{itemize}
\item[{\rm (a)}] $\mathfrak{F}$ is a birational map of $V$ onto its image {\rm ;}
\item[{\rm (b)}] $\mathfrak{F}$ has non-degenerate Jacobian dual rank and the latter equals
the embedding dimension of the affine cone over $V\,${\rm ;}
\item[{\rm (c)}] For every irreducible component $W$ of $V$, the restriction of
$\mathfrak{F}$ to $W$ has non-degenerate Jacobian dual rank and the latter equals
the embedding dimension of the affine cone over $W$.
\end{itemize}
Moreover, when these conditions hold then

\medskip

{\rm (i)} 
The inverse map to $\mathfrak{F}$ is given by the coordinates of any sufficiently general homogeneous vector of positive degree in the null
space of a weak Jacobian dual matrix of $\mathfrak{F}$.

\medskip

{\rm (ii)}
 If besides, $V$ is irreducible, the coordinates of the inverse to $\mathfrak{F}$ are
the {\rm (}ordered, signed{\rm )} $n$-minors
of an arbitrary $n\times (n+1)$ submatrix rank $n$ of a weak Jacobian dual matrix of $\mathfrak{F}$.
\end{Theorem}
(``Sufficiently general'' in the above statement means that the general $k$-linear combination of the coordinates of the vector ought
to be a non-zero-divisor on the homogeneous coordinate ring of the image.)

\begin{Remark}\rm
(1) A tall order on the computational side of the theory is to get hold of the spanning
biforms of ${\cal J}_{1,*}$ and ${\cal J}_{\, 0,*}$ without actually computing
all the generators of ${\cal J}_{\ff}$.
For the second bihomogeneous piece this is beckoningly the main problem of implicitization
in elimination theory. One of the common techniques is to try and derive the implicit equations
by means of the so-called Sylvester method starting out from the defining equations of the
symmetric algebra of $I=(\ff)$ (i.e., the bihomogeneous piece of bidegrees $(*,1)$).
In special situations, this method often throws some light on obtaining ${\cal J}_{1,*}$
-- see, e.g., \cite[Example 5.7]{syl2}.

(2) Note that the theorem gives no hint on the minimal possible degree of a representative of the inverse map.
Even in the case of a Cremona map, the degree will be elusive as the representative given by the minors will typically carry
a non-trivial factor. Thus, the question arises as to how one goes about computing the degree of the gcd of the minors
in the latter case.
Often this can be computed by a careful scrutiny of the Hilbert series (see, e.g.
\cite[Proposition 3.1]{cremona}) or by local structure (see \cite[Examples 2.4, 2.5, ]{cremona}
and also \cite[Theorem 3.1]{CremonaMexico}).

(3) It would be curious to understand the meaning of a weak Jacobian dual matrix -- or at least of its rank -- in purely geometric terms. In particular, the Jacobian dual rank seems to be related
to the field degree in a sort of complementariness.
\end{Remark}

\section{Tidbits on birational maps}

\subsection{Algebraic preliminaries}

We recall a few algebraic notions which will be used in
this part.

A very basic condition on an ideal is as
follows.
\begin{Definition}\label{idealoflineartype}\rm An ideal $I\subset R$ in a ring
$R$ is said to be {\sc of linear type} if the natural $A$-algebra
homomorphism ${\cal S}_R(I) \surjects {\cal R}_R(I)$ is injective.
\end{Definition}
Here  ${\cal S}_R(I)$ and ${\cal R}_R(I)$ denote the symmetric and
the Rees algebra of $I$, respectively. The literature on this is
quite extensive (see \cite{dual}, \cite{star}, \cite{ram1} and the
book \cite{Wolmbook}).

\medskip

Now assume that  $R=k[\xx]=k[X_0,\ldots, X_n]/{\mathfrak a}$ denotes a standard
graded $k$-algebra over a field $k$.

Let $\ff=\{f_0,\ldots,f_n\}\subset R$ be forms of a fixed degree
$d\geq 1$. We will henceforth denote by $\phi=\phi_{\ff}$ a graded
matrix whose columns generate the first syzygy module of $\ff$.
Thus, one has a free graded presentation of the homogeneous ideal
$I:=(\ff)\subset R$.
\begin{equation}\label{presentation}
\oplus_s R(-d_s) \stackrel{\phi}{\lar} R^m(-d)\rar R
\end{equation}
We denote by $\phi_1=\phi_1(\ff)$ the matrix whose columns form a
set of generators of the submodule $\mbox{\rm Im}(\phi)_{d+1}$ (i.e.,
the coordinates of such vectors as elements of the standard graded
ring $R$ are linear forms). Clearly, $\phi_1$ can always be
taken as a submatrix of $\phi$, so we informally say that $\phi_1$
is the {\sc linear part\/} of the whole syzygy matrix $\phi$. Not
infrequent is the case where $\phi_1=0$. On the other end of the
spectrum, an important special case is when $\phi_1=\phi$ - we then
say that $I$ has a {\sc linear presentation}.

A weaker condition is a certain rank requirement.
Recall that $\phi$ has a well-defined  rank $r$ if and only if $I_s(\phi)=\{0\}$
for $s>r$ and $I_r(\phi)$ has a regular element.
This is equivalent to requiring that the $R$-submodule $Z:=\mbox{\rm Im}(\phi)\subset R^m$
have rank $r$ -- i.e., that $Z_p$ be free of rank $r$ for every associated prime of $R$.
One also has that $Z$ (hence, $\phi$ too) has a rank if and only if $I$ has a regular element.

We will say that $\phi_1$ has maximal (possible) rank if $\phi$ has a rank and $\rk(\phi_1)=\rk(\phi)$.
We could also say in this situation that $\phi$
has {\sc maximal linear rank}.

\subsection{Consequences of the main criterion}

In this part we draw some corollaries of  Theorem~\ref{main_criterion_algebraic} in the case where
$\mathfrak{F}:\pp^n\dasharrow \pp^n$.
Recall that a birational map of $\pp^n$ onto itself is known as a {\sc Cremona transformation}.

We emphasize that in this context  much of the theory, including the  main criterion, gets simplified and
more visible. Cremona transformations have a uniquely defined minimal representative (both the
map and its inverse) and the corresponding degrees are important invariants in the classical theory.

Some of the applications in this section envisage solving the questions posed in \cite{bir2003}, in that
they are characteristic-free.

For example, one has the following characteristic-free
version of \cite[Theorem 4.1]{bir2003} and \cite[Proposition 4.5, (iii)]{CRS}.

\begin{Theorem}\label{linearmaximalrank}
Let $\mathfrak{F}\colon\pp^n\dasharrow \pp^m$ be a rational map, given by $m+1$
forms $\ff=\{f_0,\ldots,f_n\}$ of a fixed degree. If the image of $\mathfrak{F}$ has dimension $n$
and $\rk \phi_1=n$ {\rm (}maximal possible{\rm )} then $\mathfrak{F}$ is birational onto its image.
\end{Theorem}
\demo Following the basic {\em leitmotif\/} of \cite[Proposition 1.1]{dual},
write $\YY\cdot \phi_1=\XX\cdot \rho$,
where $\rho$ is a uniquely defined matrix
whose entries are linear forms in $k[\YY]=k[Y_0,\ldots,Y_m]$.
Note that the left (respectively, right) side
generates the presentation ideal of the symmetric algebra ${\cal S}_{k[\XX]}({\rm coker}(\phi_1))$
(respectively, ${\cal S}_{k[\YY]}({\rm coker}(\rho))$).
Therefore, the identity map
of $k[\XX,\YY]$ induces a $k$-isomorphism between these algebras.

Note the surjection ${\rm coker}(\phi_1)\surjects I=(\ff)$, whose kernel
is torsion since  $\phi_1$ has rank $n$ by hypothesis.
Therefore
\begin{equation}\label{rees_cokernel}
{\mathcal R}_{k[\XX]}({\rm coker}(\phi_1))\simeq {\mathcal R}_{k[\XX]}(I).
\end{equation}
Let $\mathcal{J}$ denote the presentation ideal of the Rees algebra ${\mathcal R}_{k[\XX]}(I)$ on $k[\XX,\YY]$.
Write $S:=k[\YY]/{\mathfrak b}$  for the homogeneous coordinate ring of the image of $\mathfrak{F}$
on $\pp^m$ and. Recall from previous sections that ${\mathfrak b}\subset \mathcal{J}$.

Write $E:={\rm coker}(\rho)$.
Since the symmetric algebra is functorial, one has
$${\cal S}_{S}(E\otimes_{k[\YY]} S)\simeq {\cal S}_{k[\YY]}(E)/{\mathfrak b}\,{\cal S}_{k[\YY]}(E)\simeq
k[\XX,\YY]/({\mathfrak b}, I_1(\XX\cdot \rho))= k[\XX,\YY]/({\mathfrak b}, I_1(\YY\cdot \phi_1)).$$
Thus, there is a surjection ${\cal S}_{S}(E\otimes_{k[\YY]} S)\surjects k[\XX,\YY]/\mathcal{J}$.

{\bf Claim}: the surjection ${\cal S}_{S}(E\otimes_{k[\YY]} S)\surjects k[\XX,\YY]/\mathcal{J}$
induces an isomorphism
$$\mathcal{R}_S(E\otimes_{k[\YY]} S)\simeq k[\XX,\YY]/\mathcal{J}.$$

Recalling that the Rees algebra $\mathcal{R}_S(E\otimes_{k[\YY]} S)$ is obtained thereof
by killing the $S$-torsion, we write
$\mathcal{T}$ for the lifting of the latter to $k[\XX,\YY]$.
Thus
$$\mathcal{R}_S(E\otimes_{k[\YY]} S)\simeq k[\XX,\YY]/({\mathfrak b}, I_1(\YY\cdot \phi_1),\mathcal{T}).$$
Then let $G\in \mathcal{T}$. Then $F(\YY)\,G\in ({\mathfrak b}, I_1(\YY\cdot \phi_1)\subset \mathcal{J}$
for some $F(\YY)\in k[\YY]\setminus {\mathfrak b}$.
Since $ \mathcal{J}$ is a prime ideal and $F(\YY)\not\in  \mathcal{J}$ because $k[\YY]\cap  \mathcal{J}={\mathfrak b}$, necessarily $G\in  \mathcal{J}$.
This proves the inclusion $({\mathfrak b}, I_1(\YY\cdot \phi_1),\mathcal{T})\subset \mathcal{J}$ as
$({\mathfrak b}, I_1(\YY\cdot \phi_1))\subset \mathcal{J}$.

Conversely, let $G=G(\XX,\YY)\in \mathcal{J}$.
By (\ref{rees_cokernel}), $\mathcal{J}$ is the $k[\XX]$-torsion of $k[\XX,\YY]/I_1(\YY\cdot\phi_1)$
lifted to $k[\XX,\YY]$.
Therefore, there exists $F(\XX)\in k[\XX]\setminus I_1(\YY\cdot\phi_1)$ -- i.e., $F(\XX)\neq 0$ --
such that $F(\XX)G(\XX,\YY)\subset I_1(\YY\cdot\phi_1)$.
Suppose that $G\not\in ({\mathfrak b}, I_1(\YY\cdot \phi_1),\mathcal{T})$.
The latter is a prime ideal, hence $F(\XX)\in ({\mathfrak b}, I_1(\YY\cdot \phi_1),\mathcal{T})$.
By the definition of $\mathcal{T}$ there exists $H(\YY)\in k[\YY]\setminus {\mathfrak b}$
such that $H(\YY)F(\XX)\in ({\mathfrak b}, I_1(\YY\cdot \phi_1))$.
Evaluating $\YY\mapsto \ff$ would give $H(\ff)F(\XX)= 0$ whence $F(\XX)=0$
since $H(\ff)\neq 0$; this is a contradiction.

To conclude we count dimensions according to \cite{ram1}, thereby finding
$$\dim S+\rk_S(E\otimes_{k[\YY]} S)= n+1 + 1.$$
Since $\dim S=n+1$ by assumption, we get $\rk_S(E\otimes_{k[\YY]} S)= 1$, hence $\rk_S(\rho\otimes_{k[\YY]} S)=n$.
Since Im$(\rho\otimes_{k[\YY]} S)\subset {\rm Im}(\psi)$ for a weak Jacobian dual matrix $\psi$ over $S$,
Theorem~\ref{main_criterion_algebraic} applies.
\qed

\medskip

Easy examples show that even if $m=n$ and the forms $\ff$ generate a linearly presented
ideal -- i.e.,  $\phi$ is linear -- they can turn out to be algebraically dependent. Perhaps the
simplest is the map $\mathfrak{F}\colon\pp^3 \dasharrow\pp^3$ defined by the
$2$-forms $X_0X_2, X_0X_3, X_1X_2, X_1X_3$.
At the other end of the spectrum, it is obvious (and easy) that
algebraic independence does not imply anything in particular concerning the rank
of $\phi_1$.
Therefore the two conditions of the proposition are pretty much mutually
independent.

The following example gives an illustration of the above proposition.

\begin{Example}\rm Consider the map $\mathfrak{F}\colon\pp^5 \dasharrow\pp^5$ defined by the
$3$-forms
$$X_0X_1X_2,\, X_0X_2X_3,\, X_0X_4X_5,\, X_1X_2X_4,\, X_2X_3X_5,\,
X_3X_4X_5.$$ Since these are square-free monomials, it is fairly
easy to compute the matrix $\phi_1$, which is here of maximal rank
($=5$). In fact, $\phi_1$ generates a free submodule. A computation
with {\it Macaulay\/} (\cite{Macaulay}) yields that (notation as
above) the three algebras ${\cal S}({\rm coker}(\phi_1)), {\cal S}(I)$ and ${\cal
R}(I)$ are all different though of the same dimension. In particular,
$I$ is not of linear type. The base ideal of the
inverse of $\mathfrak{F}$ is generated by $4$-forms, so this degree is one less
than the expected degree were $I$ of linear type (cf.
Theorem~\ref{lineartype} and its supplement).
\end{Example}
In the case of rational maps defined by monomials as this, there is a blind (i.e., algorithmic)
way to both testing birationality and obtaining the coordinates of the inverse map (see
\cite{CremonaMexico}).

In the remaining of the section, following the classical terminology,
we say that a Cremona map has degree $r$
if it is defined by forms of degree $r$ generating an ideal of codimension $\geq 2$.

One can  prove the converse to the linear
obstruction result (\cite[Proposition 3.5]{bir2003}) in
any characteristic.

\begin{Proposition}\label{lineartype}
Let $\mathfrak{F}\colon\pp^n\dasharrow \pp^n$ be a rational map whose base
ideal $I$ is of linear type. The following conditions are
equivalent:
\begin{enumerate}
\item[{\rm (a)}] $\mathfrak{F}$ is a Cremona transformation
\item[{\rm (b)}] The linear syzygy matrix $\phi_1$ has rank $n$.
\end{enumerate}
Moreover, if these conditions take place,  the inverse
map $\mathfrak{F}^{-1}$ has degree $\leq n$.
\end{Proposition}
\demo The implication (a) $\Rightarrow$ (b) has already been
observed in \cite[Proposition 3.5]{bir2003}.

The reverse implication (b) $\Rightarrow$  (a) comes from Theorem~\ref{linearmaximalrank} since the linear type
assumption for an
ideal generated by forms of the same degree implies at least the
algebraic independence of these forms.

The last assertion follows from the corresponding assertion in
Theorem~\ref{main_criterion_algebraic}.
\qed

\bigskip

The conjectured statements \cite[Corollaries 3.8 through
3.10]{bir2003} become true statements in any characteristic which we briefly restate for
convenience and reference.

\begin{Corollary}\label{onedimensional}
Let $\mathfrak{F}\colon\pp^n\dasharrow \pp^n$ be a rational map whose base
ideal $I$ is a one-dimensional locally complete
intersection at its minimal primes. Then $\mathfrak{F}$ defines a Cremona
transformation of $\pp^n$ if and only if the linear syzygy matrix $\phi_1$
has rank $n$.
\end{Corollary}
\demo As in [loc.cit.], one only needs to recall that the
assumption in the one-dimensional case means that the ideal is of
linear type (see \cite[Proposition 3.7]{conormal}, also \cite[Remark
10.5]{Trento}).
\qed

\medskip

In particular, this applies to the case of plane rational
maps, thereby yielding a strong outcome:

\begin{Corollary}\label{plane_maps}
Let $\mathfrak{F}\colon\pp^2\dasharrow \pp^2$ be a rational map
whose base ideal scheme is a codimension $2$ locally complete
intersection.
Then $\mathfrak{F}$ is a Cremona map if and only if is defined by the $2\times 2$
minors of a $3\times 2$ matrix of linear forms.
\end{Corollary}
\demo
Suppose that $\mathfrak{F}$ is Cremona. By Proposition~\ref{lineartype},
the inverse map is of degree $\leq 2$.
Since we are in the plane case, $\mathfrak{F}$ has also degree $\leq 2$.
Therefore, as the base ideal is of codimension $2$, the map is
one of the three classical types of quadratic Cremona maps,
in which case the result is directly verified.

The converse follows directly from the previous corollary.
\qed

\begin{Question}\label{dolg}
Let $f\in k[X_0,X_1,X_2]$ be a reduced homogeneous polynomial whose polar map
is a Cremona map. Is the ideal generated by the partial derivatives
of $f$ locally a complete intersection at its minimal primes?
\end{Question}
The answer is positive in characteristic zero due to the theorem
of Dolgachev (\cite{Dol}).
However, even in null  characteristic an independent proof would be of
interest as it would imply Dolgachev's result by way of Corollary~\ref{plane_maps}.

\begin{Remark}\label{entries_are_primary}\rm In the plane case ($n=2$), the
above local condition on the base ideal  is also equivalent
to requiring that its $n$th Fitting ideal $I_1(\phi)$
have codimension three, so the condition is easily verified
by computation.
\end{Remark}

\subsection{The semilinear generation degree}

In this part we raise a few  questions related to certain special
behavior of the rational map vis-\`a-vis its base ideal $I\subset
R$.

\begin{Proposition}\label{semilinearmaximalheight}
Let $\mathfrak{F}:\pp^n\dasharrow \pp^m$ be a rational map and let $W\subset \pp^m$ denote
its image.
If the codimension of the ideal $({\cal J}_{1,*})$ over the ring $(k[\YY]/I(W))[\XX]$
is at least $n$ then $\mathfrak{F}$  is birational onto $W$.
\end{Proposition}
\demo By Theorem~\ref{main_criterion_geometric}, it suffices to show that the rank of the Jacobian dual matrix $\psi$
over $S:=k[\YY]/I(W)$ is at least $n$ (hence $=n$).

For this note that if $C:={\rm coker}_S(\psi^t)$, then the symmetric algebra
${\cal S}_S(C)\simeq S[\XX]/I_1(\XX\cdot\psi^t)$.
In addition, by the very definition, $({\cal J}_{1,*})S[\XX]=I_1(\XX\cdot\psi^t)$.
It follows that
$$\dim {\cal S}_S(C)\leq \dim S+n+1-\hht ({\cal J}_{1,*})S[\XX]\leq \dim S+n+1-n=\dim S +1.
$$
On the other hand, $\dim {\cal S}_S(C)\geq \dim {\cal R}_S(C)=\dim S+\rk _S(C)$,
where ${\cal R}_S(C)$ stands for the Rees algebra of $C$.
From this $\rk_S(C)\leq 1$, hence $\rk _S(\psi)=\rk _S(\psi^t)\geq n$.
 \qed

 \bigskip

Regarding the above proof, one sees that since $\rk _S(\psi)=n$ then actually $\rk _S(C)= 1$.
In this notation, the inverse map to $\mathfrak{F}$ is defined by the coordinates of any nonzero homogenous
vector in the dual module $C^*$.

The converse statement to the effect that $\rk _S(\psi)=n$ over $k[\YY]/I(W)$
implies that the codimension of $({\cal J}_{1,*})$ over $(k[\YY]/I(W))[\XX]$
is at least $n$ is generally false, even for Cremona maps.
The reason is simply given by the fact that there are Cremona maps for which the symmetric algebra
of the base ideal has larger dimension than its Rees algebra.
No such examples exist for plane Cremona maps as the base ideal is an almost complete intersection.
For space Cremona maps, however, even monomial involutions may show this behavior as the following
example conveys.

\begin{Example}\rm
Let $\mathfrak{F}=(X^4:X^2YW:XY^2Z:Y^3Z):\pp^3\dasharrow \pp^3$.
Then $({\cal J}_{1,*})\subset k[\XX,\YY]$ has codimension $2$, i.e., $\dim {\cal S}(I)=6$,
while $\rk \psi=3$ as the map is Cremona.
Here of course $\rk(\phi_1)<3$ (actually $\rk(\phi_1)=1$ by calculation).
We note that this Cremona map is an involution up to permutation in the source and in the target.
A different approach to verify birationality can be applied in this example since the base ideal
is generated by monomials (see \cite{CremonaMexico}).
\end{Example}

In connection to Proposition~\ref{semilinearmaximalheight}, the
following numerical invariant seems worth looking at.

\begin{Definition}\rm The {\sc semilinear generation degree} of an ideal
$I\subset R=k[\XX]$, denoted sgd$(I)$, is the minimum $s\geq 1$ such that
$({\cal J}_{1,*})=({\cal J}_{1, \,\leq s})$,
where the ideal on the right side is the subideal
generated by all (``semilinear'') relations of bidegree $(1,r)$ with
$r\leq s$.
\end{Definition}
A souped-up version would require only $\hht
({\cal J}_{1,*})=\hht ({\cal J}_{1, \,\leq s})$, in which case
the invariant could be dubbed the {\sc semilinear codimension degree}.

This numerical invariant expresses the measure of how far off through the
successive powers $I^m$ of the ideal $I$ one has to go in order to
pick a sufficient number of linear syzygies of $I^m$ to generate the
ideal $({\cal J}_{1,*})$. In case $I$ is actually the base ideal of
a Cremona map then this number may have bearing to the degree of the
coordinates of the inverse map.
Here is a very crude lower bound:

\begin{Proposition}\label{crude_sgd}
Let $I\subset k[\XX]$ denote the base ideal of a Cremona map $\mathfrak{F}:\pp^n\dasharrow \pp^n$.
Then sgd$(I)\geq \deg(\mathfrak{F}^{-1})/n$, where $\deg \mathfrak{F}^{-1}$ denotes the
degree of the defining coordinates of the inverse map.
\end{Proposition}
\demo
Using item (ii) of Theorem~\ref{main_criterion_geometric} one immediately
has  $\deg \mathfrak{F}^{-1}\leq {\rm sgd}(I)\,n$, as required.
\qed

\medskip

For instance, if $I$ is of linear type then $\deg \mathfrak{F}^{-1}\leq n$,
a property already observed in Proposition~\ref{lineartype}.

\medskip

\begin{Question}\label{semilinearnumber}\rm
Find upper bounds for sgd$I$ in
terms of other invariants of $I$.
\end{Question}

A related question was raised by Vasconcelos in
connection with elimination theory. Note that the above numerical
invariant makes sense even for a general rational map with target
some $\pp^m$.

\begin{Question}\label{semilinearnumber_vs_eliminationdegree}\rm
Let more generally $\mathfrak{F}\colon\pp^n\dasharrow \pp^m$ be a rational map,
given by $m+1$ $k$-linearly independent forms $\ff=\{f_0,\ldots,f_m\}\subset k[\XX]$ of a fixed degree
$\geq 2$, without proper common factor. Let $I=(f_0,\ldots,f_m)$
denote as above the base ideal of $\mathfrak{F}$ and let $d$ denote the highest
degree of a minimal form defining the image of $\mathfrak{F}$ in $\pp^m$. Is sgd$(I)\leq d$?
\end{Question}

Note that the answer is trivially affirmative for ideals of fiber type.
There is both theoretical and computational evidence for an affirmative answer
in the case of an $\XX$-primary almost complete intersection (see \cite{syl2}).


\bigskip

\noindent{\bf Addresses:}

\medskip

{\sc A. V. D\'oria}, Departamento de Matem\'atica, Universidade Federal de Sergipe,
Sergipe, Brazil

{\it Email}: {\sc avsdoria@dmat.ufpe.br}

\bigskip

{\sc S. H. Hassanzadeh}, Faculty of Mathematical Sciences and Computer, Tarbiat
Moallem University, Tehran, Iran {\em and\/} Departamento de Matem\'atica, CCEN, Universidade
Federal de Pernambuco,  Pernambuco,
Brazil

{\it Email}: {\sc hamid@dmat.ufpe.br}

\bigskip

{\sc A. Simis}, Departamento de Matem\'atica, CCEN, Universidade
Federal de Pernambuco,  Pernambuco,
Brazil

{\it Email}: {\sc aron@dmat.ufpe.br}


\begin{thebibliography}{99}


\bibitem{Macaulay}{D. Bayer and M. Stillman, {\sc Macaulay}: a
computer algebra system for algebraic geometry, Macaulay version 3.0
1994 (Macaulay for Windows by Bernd Johannes Wuebben, 1996).}


\bibitem{CRS}{C. Ciliberto, F. Russo and A. Simis, Homaloidal hypersurfaces and
hypersurfaces with vanishing Hessian, Advances in Math., {\bf 218} (2008) 1759--1805.}

\bibitem{Dol}{I. V. Dolgachev, Polar Cremona transformations, Mich. Math. J.,
{\bf 48} (2000), 191--202.}


\bibitem{Trento}{J. Herzog, A. Simis and W. V. Vasconcelos, Koszul homology
and blowing-up rings, in {\sc Commutative Algebra} (S. Greco and G.
Valla, eds.),  Lecture Notes in Pure and Applied Math.  {\bf 84},
Marcel-Dekker, New York, 1983, 79--169.}



\bibitem{syl2}{J. Hong, A. Simis and W. V. Vasconcelos, The equations of
almost complete intersections, arXiv:0906.1591v1 [math.AC] 8 Jun 2009.}


\bibitem{HKS}{K. Hulek, S. Katz, F.-O. Schreyer,  Cremona transformations and
syzygies, Math. Z. {\bf 209} (1992),  419--443.}

\bibitem{PanRusso}{I. Pan and F. Russo, Cremona transformations and
special double structures, Manuscripta Math. {\bf 117} (2005), 491--510.}



\bibitem{cremona}{F. Russo and A. Simis, On birational maps and Jacobian matrices,
Compositio Math. {\bf 126} (2001), 335--358.}


\bibitem{bir2003}{A. Simis, Cremona transformations and some related
algebras, J. Algebra {\bf 280 (1)} (2004), 162--179.}

\bibitem{star}{A. Simis, B. Ulrich  and W. Vasconcelos,
Tangent star cones,  J.  reine angew. Math., {\bf 483} (1997), 23--59.}

\bibitem{dual}{A. Simis, B. Ulrich  and W. Vasconcelos, Jacobian dual
fibrations, Amer. J. Math. {\bf 115} (1993), 47--75.}

\bibitem{ram1}{A. Simis,  B. Ulrich and W. V. Vasconcelos, Rees algebras of modules,
Proc. London Math. Soc., {\bf 87 (3)} (2003), 610--646.}

\bibitem{ram2}{A. Simis,  B. Ulrich and W. V. Vasconcelos, Codimension, multiplicity
and integral extensions, Math. Proc. Camb. Phil. Soc. {\bf 130},
(2001), 237--257.}

\bibitem{conormal}{A. Simis  and W. Vasconcelos, The syzygies of the conormal module,
Amer. J. Math., {\bf 103 (2)} (1981), 203--224.}



\bibitem{SiVi} {A. Simis and
R. H.  Villarreal, Constraints for the normality of monomial subrings and
birationality, Proc. Amer. Math. Soc. {\bf 131} (2003), 2043--2048.}

\bibitem{SimisVilla}{A. Simis and R. Villarreal, Linear syzygies and birational
combinatorics, Results Math. {\bf 48} (2005), 326--343.}

\bibitem{CremonaMexico}{A. Simis and R. Villarreal, Combinatorics of monomial Cremona maps,
arXiv:0904.4065v1 [math.AG].}

\bibitem{Wolmbook}{W. V. Vasconcelos, {\sc Arithmetic of Blowup Algebras}.
London Math. Soc.,
Lecture Notes Series {\bf 195}, Cambridge University Press, Cambridge, 1994.}



\end{thebibliography}
\end{document}